\documentclass[11pt,leqno]{amsart}
\usepackage{amssymb,verbatim,enumerate,ifthen,times}
\usepackage[mathscr]{eucal}
\usepackage[utf8]{inputenc}
\usepackage[T1]{fontenc}
\usepackage{tabu}
\usepackage{dsfont}
\usepackage{accents}

\usepackage{color}
\usepackage{hyperref}

\makeatletter
\@namedef{subjclassname@2020}{\textup{2020} Mathematics Subject Classification}
\makeatother

\def\N{\mathbb{N}}

\def\Q{\mathbb{Q}}
\def\R{\mathbb{R}}

\def\cl{\mbox{\rm cl}}
\def\diam{\mbox{\rm diam}}
\def\Ref(#1|#2){(#1\hspace{.3mm}|\hspace{.3mm}#2)}

\newtheorem{theorem}{Theorem}
\newtheorem*{theorem*}{Theorem}
\long\def\Thm#1#2{\ifthenelse{\equal{#1}{*}}{\begin{theorem*}#2\end{theorem*}}
             {\begin{theorem}\label{T#1}#2\end{theorem}}}
\newtheorem{Atheorem}{Theorem}

\def\thm#1{Theorem~\ref{T#1}}

\newtheorem{proposition}[theorem]{Proposition}
\newtheorem*{proposition*}{Proposition}
\long\def\Prp#1#2{\ifthenelse{\equal{#1}{*}}{\begin{proposition*}#2\end{proposition*}}
             {\begin{proposition}\label{P#1}#2\end{proposition}}}
\def\prp#1{Proposition~\ref{P#1}}

\newtheorem{corollary}[theorem]{Corollary}
\newtheorem*{corollary*}{Corollary}
\long\def\Cor#1#2{\ifthenelse{\equal{#1}{*}}{\begin{corollary*}#2\end{corollary*}}
             {\begin{corollary}\label{C#1}#2\end{corollary}}}

\newtheorem{lemma}[theorem]{Lemma}
\newtheorem*{lemma*}{Lemma}
\long\def\Lem#1#2{\ifthenelse{\equal{#1}{*}}{\begin{lemma*}#2\end{lemma*}}
             {\begin{lemma}\label{L#1}#2\end{lemma}}}
\def\lem#1{Lemma~\ref{L#1}}

\theoremstyle{definition}
\newtheorem{definition}[theorem]{Definition}
\newtheorem*{definition*}{Definition}
\long\def\Defn#1#2{\ifthenelse{\equal{#1}{*}}{\begin{definition*}\rm #2\end{definition*}}
             {\begin{definition}\label{D#1}\rm #2\end{definition}}}

\newtheorem{remark}{Remark}
\newtheorem*{remark*}{Remark}
\long\def\Rem#1#2{\ifthenelse{\equal{#1}{*}}{\begin{remark*}\rm #2\end{remark*}}
             {\begin{remark}\label{R#1}\rm #2\end{remark}}}
\def\rem#1{Remark~\ref{R#1}}

\newtheorem{example}{Example}
\newtheorem*{example*}{Example}
\long\def\Exa#1#2{\ifthenelse{\equal{#1}{*}}{\begin{example*}\rm #2\end{example*}}
             {\begin{example}\label{Ex#1}\rm #2\end{example}}}

\def\eq#1{{\rm(\ref{E#1})}}

\def\Eq#1#2{\ifthenelse{\equal{#1}{*}}
  {\begin{equation*}\begin{aligned}[]#2\end{aligned}\end{equation*}}
  {\begin{equation}\begin{aligned}\label{E#1}#2\end{aligned}\end{equation}}}

\begin{document}
\date{\today}

\title{A Pexider equation containing the aritmetic mean}
\author[T. Kiss]{Tibor Kiss}
\address{
Institute of Mathematics, University of Debrecen,\newline
\indent P. O. Box 400, H-4002 Debrecen, Hungary \newline
\indent ELKH-DE Equations, Functions, Curves and their Applications
Research Group}
\email{kiss.tibor@science.unideb.hu}
 \subjclass[2020]{Primary 39B22, Secondary 26E60}
\keywords{Functional equation, Pexider equation, arithmetic mean.}
\thanks{The research of the author was supported in part by the ÚNKP-22-4 New National 
Excellence Program of the Ministry for Culture and Innovation from the source of the National 
Research, Development and Innovation Fund, in part by NKFIH Grant K-134191, and in part by the Eötvös Loránd Research Network (ELKH)}

\begin{abstract}
In this paper we determine the solutions $(\varphi,f_1,f_2)$ of the Pexider functional equation
\Eq{*}{
\varphi\Big(\frac{x+y}2\Big)\big(f_1(x)-f_2(y)\big)=0,\qquad (x,y)\in I_1\times I_2,
}
where $I_1$ and $I_2$ are nonempty open subintervals. Special cases of the above equation 
regularly arise in problems with two-variable means. We show that, under a rather weak 
regularity condition, the coordinate-functions of a typical solution of the equation are 
constant over several subintervals of their domain. The regularity condition in question will 
be that the set of zeros of $\varphi$ is closed. We also discuss particular solutions where 
this condition is not met.

\end{abstract}

\maketitle

\section{Introduction}
Functional equations of type
\Eq{E0}{
\varphi\Big(\frac{x+y}2\Big)\big(f_1(x)-f_2(y)\big)=0,\qquad (x,y)\in I_1\times I_2,
}
are closely connected with various problems with two-variable means. In relation to the 
Matkowski--Sutô problem \cite{DarPal02}, equation \eq{E0} was 
studied assuming that $I_1=I_2$, $g:=f_1=f_2$, $\varphi=g+\lambda$ for some constant 
$\lambda\in\R$, and that $g$ is of the form $d\circ c$, where $d$ is a derivative and $c$ is 
continuous.

In paper \cite{DarTot15} Z. Daróczy and V. Totik investigated a co-equation of 
\eq{E0} and proved that continuous differentiability of the solutions implies that they are 
constant. The authors also showed that there exist differentiable, nowhere monotone (and hence 
non-constant) solutions as well. The interested reader can also find equations with a similar 
spirit in \cite{DarLac08}.

In \cite{KisPal19}, as an auxiliary functional equation, a special case of \eq{E0} played 
an important role in solving the equality problem of Cauchy means and quasi-arithmetic means. 
The solutions of \eq{E0} were determined assuming that $I_1=I_2$, $f_1=f_2$, and that the set 
of zeros of $\varphi$ is closed. Also in terms of functions and domains, equation \eq{E0} is a 
pexiderization of the auxiliary equation of \cite{KisPal19}.

Note that, in fact, for given functions $f_1$ and $f_2$, it is not so difficult to construct 
$\varphi$ to  obtain a solution. Indeed, let $f_1:I_1\to\R$ and $f_2:I_2\to\R$ be any functions, 
$\lambda\in f_1(I_1)$ be any value but fixed, and define $\Lambda:=f_2^{-1}(\{\lambda\})$. Then  
$(\varphi,f_1,f_2)$ is a solution of \eq{E0} provided that $\varphi(u)=0$ whenever  
$u\in\frac12\big(I_1+(I_2\setminus\Lambda)\big)$. The aim of this paper is to solve the equation 
under a weak but reasonable condition that guarantees that the solutions have regular, more 
precisely, constant parts. We are going to describe the solution set of functional equation 
\eq{E0}, where $I_1,I_2\subseteq\R$ stand for nonempty open intervals, the functions 
$\varphi:D:=\frac12(I_1+I_2)\to\R$, $f_1:I_1\to\R$, and $f_2:I_2\to\R$ are unknown, and the set 
of zeros of $\varphi$ is assumed to be closed in $D$. Note that this latter condition is 
trivially satisfied whenever $\varphi$ is continuous, strictly monotone, injective or, say, has  finitely many zeros.

Adapted from paper \cite{KisPal19}, for $H\subseteq\R$ and a given function $\ell:H\to\R$, the 
set $\{u\in H\mid\ell(u)=0\}$ will be denoted by $\mathscr{Z}_\ell$. For its complementary set 
$H\setminus\mathscr{Z}_\ell$ we shall use $\mathscr{Z}_\ell^c$. Whenever we write 
$(\varphi,f_1,f_2)$, we mean that $\varphi$, $f_1$, and $f_2$ are real valued functions defined 
on $D$, $I_1$, and $I_2$, respectively.  A triplet of functions $(\varphi,f_1,f_2)$ will be 
called \emph{extremal} if $D\subseteq\mathscr{Z}_\varphi$, that is, if $\varphi$ is identically 
zero or there exists $\lambda\in\R$ such that $f_1(x)=f_2(y)=\lambda$ for all $x\in I_1$ and 
$y\in I_2$. We can then make the following trivial observation.

\Prp{nec}{If $(\varphi,f_1,f_2)$ is extremal then $(\varphi,f_1,f_2)$ solves equation \eq{E0}.}

In fact, it is the above statement that motivates the name extremal. As we will see later, a 
general solution of \eq{E0}, roughly speaking, is a mixture of extremal solutions, and this nice 
property is provided by the condition that $\mathscr{Z}_\varphi$ is closed.

A nonempty family $\{U\subseteq\R\mid U\text{ is an interval}\}$ will be called a \emph{sequence 
of intervals}, if it is countable. Then the family in question will be simply referred to as 
$(U_n)$, that is, we do not indicate the index set.

For given subsets $P,Q,S\subseteq\R$, define
\Eq{*}{
(S \hspace{.3mm}|\hspace{.3mm} P)_Q:=(2P-S)\cap Q.
}
In words, $(S \hspace{.3mm}|\hspace{.3mm} P)_Q$ consists of those elements of $Q$ that are 
reflections of some element 
of $S$ with respect some element of $P$. If at least one of the sets is empty, 
then $(S \hspace{.3mm}|\hspace{.3mm} P)_Q$ is empty. For brevity, when reflecting to a given point, that 
is, if $P$ is of the form $\{p\}$, then instead of $(S \hspace{.3mm}|\hspace{.3mm} \{p\})_Q$ we simply write 
$(S \hspace{.3mm}|\hspace{.3mm} p)_Q$. Finally, having an interval $J\subseteq\R$, we shall say 
that $H\subseteq J$ is a \emph{proper} subset if $H\notin\{\emptyset,J\}$.

Later we want to refer to some trivial properties of a set of the form $(S 
\hspace{.3mm}|\hspace{.3mm} P)_Q$, hence we formulate the following lemma.

\Lem{L0}{Let $P,Q,S\subseteq\R$ be any subsets.
\begin{enumerate}[(i)]\itemsep=1mm
\item If $P_1\subseteq P_2\subseteq\R$, $Q_1\subseteq Q_2\subseteq\R$, and $S_1\subseteq 
S_2\subseteq\R$, then
\Eq{*}{
(S_1 \hspace{.3mm}|\hspace{.3mm} P_1)_{Q_1}\subseteq(S_2 \hspace{.3mm}|\hspace{.3mm} 
P_2)_{Q_2}.
}
\item The sets $(S \hspace{.3mm}|\hspace{.3mm} P)_Q$ and $(Q \hspace{.3mm}|\hspace{.3mm} P)_S$ are mirror images of each other with 
respect to $P$, more precisely, we have
\Eq{*}{
((S \hspace{.3mm}|\hspace{.3mm} P)_Q \hspace{.3mm}|\hspace{.3mm} P)_S=(Q \hspace{.3mm}|\hspace{.3mm} P)_S.
}
\item If $Q$ is open and at least one of the sets $P$ and $S$ is open, then $(S \hspace{.3mm}|\hspace{.3mm} P)_Q$ 
is open.
\item Finally, if $P,Q$, and $S$ are intervals, so is $(S \hspace{.3mm}|\hspace{.3mm} P)_Q$ and
\Eq{*}{
\inf(S \hspace{.3mm}|\hspace{.3mm} P)_Q=\max(\inf Q,2\inf P-\sup S)\qquad\text{and}\qquad
\sup(S \hspace{.3mm}|\hspace{.3mm} P)_Q=\min(\sup Q,2\sup P-\inf S),
}
provided that $(S \hspace{.3mm}|\hspace{.3mm} P)_Q$ is nonempty.
\end{enumerate}}

\begin{proof}
Assertions (i), (iii), and (iv) are direct consequences of the definition.

The proof of (ii) is also elementary. Set $A:=((S \hspace{.3mm}|\hspace{.3mm} P)_Q 
\hspace{.3mm}|\hspace{.3mm} P)_S$ and 
$B:=(Q \hspace{.3mm}|\hspace{.3mm} P)_S$. By  definition, $(S \hspace{.3mm}|\hspace{.3mm} 
P)_Q\subseteq Q$, hence, in view of statement (i), $A\subseteq B$ follows. The reverse 
inclusion is trivial if $B$ is empty, thus we may assume that this is not the case. Let $x\in B$ 
be any point. Then $x\in S$ and there exist $y\in Q$ and $p\in P$ such that $x=2p-y$. Thus 
$2p-x=y\in(S \hspace{.3mm}|\hspace{.3mm} P)_Q$, which, together with $x\in S$, yields that $x\in 
2p-(S \hspace{.3mm}|\hspace{.3mm} P)_Q\subseteq A$.
\end{proof}

We will use what the next lemma states directly in our theorem about the solutions of \eq{E0}. 
Briefly, it says that at least one of the mirror images of $I_1$ and $I_2$ for some nonempty 
subinterval of $D$ is always located at one end of $I_2$ and $I_1$, respectively. Note that the 
subinterval for which we are reflecting may consist of a single point.

\Lem{L1}{Let $H\subseteq D$ be a nonempty subinterval and $i,j\in\{1,2\}$ with $i\neq j$. Then 
$(I_i \hspace{.3mm}|\hspace{.3mm} H)_{I_j}$ is a nonempty open subinterval of $I_j$ furthermore
\begin{enumerate}[(i)]\itemsep=1mm
\item if $\inf(I_i \hspace{.3mm}|\hspace{.3mm} H)_{I_j}\in I_j$ then $\sup(I_j 
\hspace{.3mm}|\hspace{.3mm} H)_{I_i}=\sup I_i$ and
\item if $\sup(I_i \hspace{.3mm}|\hspace{.3mm} H)_{I_j}\in I_j$ then $\inf I_i=\inf(I_j 
\hspace{.3mm}|\hspace{.3mm} H)_{I_i}$.
\end{enumerate}
}

\begin{proof} The fact that $(I_i \hspace{.3mm}|\hspace{.3mm} H)_{I_j}$ is nonempty follows from 
the inclusion $H\subseteq D$.

To prove (i), assume that $\inf(I_i \hspace{.3mm}|\hspace{.3mm} H)_{I_j}\in I_j$. Then, by (iv) 
of \lem{L0}, $\inf(I_i \hspace{.3mm}|\hspace{.3mm} H)_{I_j}=2\inf H-\sup I_i\in\R$, and hence 
$\inf H,\sup I_i\in\R$. Then, in view of (ii) and (iv) of \lem{L0}, we have
\Eq{*}{
\sup(I_j \hspace{.3mm}|\hspace{.3mm} H)_{I_i}&
=\sup((I_i \hspace{.3mm}|\hspace{.3mm} H)_{I_j} \hspace{.3mm}|\hspace{.3mm} H)_{I_i}
=\min(\sup I_i,2\sup H-\inf(I_i \hspace{.3mm}|\hspace{.3mm} H)_{I_j})\\&
=\min(\sup I_i,2\sup H-2\inf H+\sup I_i)
=\min(\sup I_i,\sup I_i+2\diam(H))=\sup I_i,}
since $\sup I_i\in\R$ and $0\leq\diam(H)\leq+\infty$. Statement (ii) is similarly verified.
\end{proof}

\section{Solutions where the interior of $\mathscr{Z}_\varphi$ is empty}

In this section we characterize the solution set of \eq{E0} in the case where the interior of 
$\mathscr{Z}_\varphi$ is empty. To do this, we need the following property of the members of 
the support of $\varphi$.

\Prp{P0}{Let $(\varphi,f_1,f_2)$ be a solution of equation \eq{E0} such that 
$\mathscr{Z}_\varphi$ is different from but closed in $D$. Then, for any 
$p\in\cl\,\mathscr{Z}_\varphi^c$, there exists a constant $\lambda\in\R$ such that
\Eq{Con}{
f_1(x)=f_2(y)=\lambda,\qquad x\in(I_2 \hspace{.3mm}|\hspace{.3mm} p)_{I_1}\text{ and }y\in(I_1 \hspace{.3mm}|\hspace{.3mm} p)_{I_2}.
}}

\begin{proof}
Let $p\in\cl\,\mathscr{Z}_\varphi^c$. Then, by \lem{L1}, 
$A:=(I_2 \hspace{.3mm}|\hspace{.3mm} p)_{I_1}$ and $B:=(I_1 \hspace{.3mm}|\hspace{.3mm} p)_{I_2}$ are nonempty open intervals in $I_1$ and 
$I_2$, respectively.

Assume first that $p\in\mathscr{Z}_\varphi^c$. We are going to show that any member of 
$A$ has an open neighbourhood in $A$ on which $f_1$ is constant. Let $x_0\in A$ be any 
point. Then there exists $y_0\in B$, such that $p=\frac12(x_0+y_0)$. Using that the sets $A$ 
and $\mathscr{Z}_\varphi^c$ are open, one can find $r>0$ such that
\Eq{*}{
U(x_0,r):=\,]x_0-r,x_0+r[\,\subseteq A
\qquad\text{and}\qquad
]p-\tfrac{r}2,p+\tfrac{r}2[\,\subseteq\mathscr{Z}_\varphi^c.
}
Applying equation \eq{E0} for $x=x_0+t$ with $|t|<r$ and for $y_0$, we get that
\Eq{*}{
\varphi\big(p+\tfrac{t}2\big)\big(f_1(x)-f_2(y_0)\big)=0,
}
where $\varphi\big(p+\tfrac{t}2\big)\neq 0$. Consequently, with $\lambda:=f_2(y_0)$, we have 
$f_1(x)=\lambda$ whenever $x\in U(x_0,r)$. That $x_0$ was arbitrary, implies that 
$f_1(x)=\lambda$ if $x\in A$. By the same reasoning we obtain that $f_2$ 
must also be constant on $B$ with some value $\mu\in\R$. Notice that $y_0\in B$, thus 
$\mu=\lambda$ follows.

Assume now that $p\in\cl\,\mathscr{Z}_\varphi^c$ with $\varphi(p)=0$ and let 
$x_0\in A$. Then $p=\frac12(x_0+y_0)$ for some $y_0\in B$. Let 
$(p_n)\subseteq\mathscr{Z}_\varphi^c$ be a sequence tending to $p$. We claim that there 
exists $n_0\in\N$ such that $x_0\in(I_2 \hspace{.3mm}|\hspace{.3mm} p_{n_0})_{I_1}$.

Indeed, otherwise, an infinite number of the intervals in question would always be to the left 
or right of the point $x_0$, more precisely, there would exist a subsequence $(p_{n_k})$ 
such 
that, for all $k\in\N$, either
\Eq{*}{
\sup(I_2 \hspace{.3mm}|\hspace{.3mm} p_{n_k})_{I_1}\leq x_0
\qquad\text{or}\qquad
x_0\leq \inf(I_2 \hspace{.3mm}|\hspace{.3mm} p_{n_k})_{I_1}.
}

The first inequality yields that $2p_{n_k}-\inf I_2\leq x_0$ for all $k\in\N$. Thus, 
particularly, $\inf I_2\in\R$ follows. Taking 
$k\to\infty$, we get that
\Eq{*}{
p\leq\tfrac12(x_0+\inf I_2)<\tfrac12(x_0+y_0)=p.
}
The second inequality leads to a similar type of contradiction. Thus 
$x_0\in(I_2 \hspace{.3mm}|\hspace{.3mm} p_{n_0})_{I_1}$ for some $n_0\in\N$.

In view of the first part of the proof, it follows that $f_1$ is constant on some open 
neighbourhood of $x_0$. But $x_0$ was arbitrary, hence $f_1$ is identically $\lambda$ on $A$ for 
some constant $\lambda\in\R$. A similar argument shows that $f_2$ must take some $\mu\in\R$ on 
$B$.

All that remains is to show is that $\lambda$ and $\mu$ are equal. Let $x_0\in A$ and $y_0\in 
B$ 
be such that $p=\frac12(x_0+y_0)$, and, again, pick a sequence 
$(p_n)\subseteq\mathscr{Z}_\varphi^c$ tending to $p$. Then there exists $n_0\in\N$ such that 
$y_0\in(A \hspace{.3mm}|\hspace{.3mm} p_{n_0})_{I_2}$. Otherwise, for some 
subsequence $(p_{n_m})$, we would have either
\Eq{*}{
\sup(A \hspace{.3mm}|\hspace{.3mm} p_{n_m})_{I_2}\leq y_0\qquad\text{or}\qquad
y_0\leq\inf(A \hspace{.3mm}|\hspace{.3mm} p_{n_m})_{I_2}.}
Both inequalities contradict the fact that $x_0\in A$. Consequently, there exist $n_0\in\N$ 
and $x\in A$ such that $p_{n_0}=\frac12(x+y_0)$. Then, by equation 
\eq{E0}, we have $\lambda=f_1(x)=f_2(y_0)=\mu$, which finishes the proof.
\end{proof}

\Thm{T0}{Let $\varphi:D\to\R$ be such that $\mathscr{Z}_\varphi$ is closed in $D$ and has an 
empty interior. Then $(\varphi,f_1,f_2)$ solves functional equation 
\eq{E0} if and only if $(\varphi,f_1,f_2)$ is extremal.}

\begin{proof}
The sufficiency is obvious, so we address the necessity.

If $\mathscr{Z}_\varphi^c$ is empty, then $D\subseteq\mathscr{Z}_\varphi$ and we are done. 
Hence assume that this is not the case.

By \prp{P0}, one can immediately obtain that any member of 
$(I_2 \hspace{.3mm}|\hspace{.3mm} \mathscr{Z}_\varphi^c)_{I_1}$ has an open neighbourhood 
in $(I_2 \hspace{.3mm}|\hspace{.3mm} \mathscr{Z}_\varphi^c)_{I_1}$ on which $f_1$ is constant. This yields that $f_1$ 
must be constant on any component of $(I_2 \hspace{.3mm}|\hspace{.3mm} 
\mathscr{Z}_\varphi^c)_{I_1}$. An analogous 
statement holds with $f_2$ and $(I_1 \hspace{.3mm}|\hspace{.3mm} \mathscr{Z}_\varphi^c)_{I_2}$.

Now observe that for any components $K\subseteq(I_2 \hspace{.3mm}|\hspace{.3mm} \mathscr{Z}_\varphi^c)_{I_1}$ and 
$L\subseteq(I_1 \hspace{.3mm}|\hspace{.3mm} \mathscr{Z}_\varphi^c)_{I_2}$, the constant 
functions $f_1|_K$ and $f_2|_L$ take a common value. Indeed, if this were not so, then, using 
equation \eq{E0}, one could deduce that $\varphi$ is identically zero on $\frac12(K+L)$, 
contradicting that $\mathscr{Z}_\varphi$ has an empty interior. Let us denote the common value 
taken by the functions $f_1|_{(I_2 \hspace{.3mm}|\hspace{.3mm} \mathscr{Z}_\varphi^c)_{I_1}}$ 
and 
$f_2|_{(I_1 \hspace{.3mm}|\hspace{.3mm} \mathscr{Z}_\varphi^c)_{I_2}}$ by $\lambda$.

Now take a point $x\in I_1$ arbitrarily and let $y\in I_2$ be such that 
$\frac12(x+y)\in\mathscr{Z}_\varphi^c$. We note that such $y$ exists: otherwise 
$\frac12(x+I_2)\subseteq\mathscr{Z}_\varphi$ would be true, which is a contradiction. Thus, by 
equation \eq{E0}, $f_1(x)=f_2(y)$ follows, where
\Eq{*}{y\in (2\mathscr{Z}_\varphi^c-x)\cap I_2\subseteq(I_1 \hspace{.3mm}|\hspace{.3mm} 
\mathscr{Z}_\varphi^c)_{I_2}.} 
Consequently, $f_1(x)=\lambda$. An analogous argument shows that the same is true for $f_2$, 
that is, $f_2(y)=\lambda$ if $y\in I_2$.
\end{proof}

\Rem{R0}{In view of \lem{L1} and \prp{P0}, for a solution of \eq{E0}, either at least one of the 
functions $f_1$ and $f_2$ is constant on its domain or at the same end of their domain they are 
constant with the same value. Furthermore if $\mathscr{Z}_\varphi$ has an empty interior or, 
equivalently, if $\mathscr{Z}_\varphi^c$ is dense in $D$, the only possible solution is the 
extremal one. To have these statements, we heavily use the topological condition on 
$\mathscr{Z}_\varphi$. The following example demonstrates that without this topological 
assumption the situation may be completely different.
	
Let $\psi:D\to\R$ be any function,
\Eq{*}{
f_1(x):=\chi_{\Q}(x),\,\,\, \text{if }x\in I_1,\quad
f_2(y):=\chi_{\Q}(y),\,\,\, \text{if }y\in I_2,
\quad\text{and}\quad
\varphi(u):=\chi_{\Q}(u)\psi(u),\,\,\, \text{if }u\in D,
}
where $\chi_{\Q}:\R\to\R$ stands for \emph{the characteristic function of rational numbers}, 
that is,
\Eq{*}{
\chi_\Q(t)=
\begin{cases}
1&\text{if }t\in\Q,\\
0&\text{if }t\in\R\setminus\Q.
\end{cases}
}
	
Then $(\varphi,f_1,f_2)$ solves equation \eq{E0} and $\mathscr{Z}_\varphi$ is closed in $D$ if and only if $(\varphi,f_1,f_2)$ is extremal or, equivalently, if $\psi$ vanishes at any rational member of $D$.
	
Indeed, let $x\in I_1$ and $y\in I_2$. If $x,y\in\Q$ or $x,y\in\R\setminus\Q$, then 
$f_1(x)=f_2(y)$. If exactly one of them is rational, then $p:=\frac12(x+y)$ must be irrational, that is $\varphi(p)=0$. The second part of the statement is trivial in view of the inclusion $D\cap(\R\setminus\Q)\subseteq\mathscr{Z}_\varphi$.
}

\section{The solution set in general}

For an index $k\in\{1,2\}$ and a constant $\lambda\in\R$, define
\Eq{*}{
a_k(\lambda):=\sup\{u\in I_k\mid f_k(t)=\lambda\text{ if }t\in\,]\inf I_k,u[\,\}
\,\text{ and }\,
b_k(\lambda):=\inf\{u\in I_k\mid f_k(t)=\lambda\text{ if }t\in\,]u,\sup I_k[\,\}.
}

In words, the extended real numbers defined above can be interpreted as follows. The set $]\inf 
I_k, a_k(\lambda)[$ is the largest open subinterval at the left end of $I_k$ on which $f_k$ is  
identically $\lambda$. We have $a_k(\lambda)=-\infty$ if and only if $I_k$ has no subinterval of 
positive length whose left endpoint is $\inf I_k$ and on which 
$f_k$ takes $\lambda$. Finally, $a_k(\lambda)=\sup I_k$ if and only if $f_k$ takes 
$\lambda$ on the whole interval $I_k$. With obvious modifications, similar statements are true 
with $b_k(\lambda)\leq+\infty$.

Observe that, whenever $(\varphi,f_1,f_2)$ is a non-extremal solution of \eq{E0}, then obviously 
\begin{enumerate}\itemsep=1mm
 \item[(C1)] neither $f_1$ nor $f_2$ is constant on its domain or
 \item[(C2)] exactly one of the functions $f_1$ and $f_2$ is constant on its domain.
\end{enumerate}
Indeed, if $(\varphi,f_1,f_2)$ is not extremal, then neither the interior of 
$\mathscr{Z}_\varphi$ nor the complementary set $\mathscr{Z}_\varphi^c$ is non-empty. So if 
$f_1$ and $f_2$ were both constant, they would have to take the same value, which contradicts 
the fact that the solution is not extremal.

\Thm{M}{Let $\varphi:D\to\R$ be such that $\mathscr{Z}_\varphi$ is closed in $D$. Then 
$(\varphi,f_1,f_2)$ solves functional equation \eq{E0} if and only if one of the following 
statements holds true.
\begin{enumerate}[(i)]\itemsep=1mm
\item The triplet $(\varphi,f_1,f_2)$ is extremal.
\item There exist constants $\lambda,\mu\in\R$ and open intervals $U_1,U_2\subseteq I_1$ and 
$V_1,V_2\subseteq I_2$ with
\Eq{*}{
\inf U_1=\inf I_1\quad\text{and}\quad\inf V_1=\inf I_2
\qquad\text{or}\qquad
\sup U_2=\sup I_1\quad\text{and}\quad\sup V_2=\sup I_2}
such that
\Eq{*}{
f_1(x)=f_2(y)=\lambda\quad\text{if}\quad (x,y)\in U_1\times V_1\qquad\text{or}\qquad
f_1(x)=f_2(y)=\mu\quad\text{if}\quad (x,y)\in U_2\times V_2}
and
\Eq{*}{
\varphi(u)=0\quad\text{whenever}\quad u\in\tfrac12(K_1+I_2)\cup\tfrac12(I_1+K_2),}
where $K_1:=I_1\setminus(U_1\cup U_2)$ and $K_2:=I_2\setminus(V_1\cup V_2)$ are closed proper 
subintervals.
\item There exists a constant $\lambda\in\R$ and, for some index $i\in\{1,2\}$, there exists a 
sequence $(U_n)$ of disjoint nonempty open subintervals of $I_i$ such that
\Eq{e2}{
f_i(x)=\lambda\quad\text{if}\quad x\in\bigcup_n U_n,\qquad 
f_j(y)=\lambda\quad\text{if}\quad y\in I_j,
\qquad\text{and}\qquad
\varphi(u)=0\quad\text{if}\quad u\in\tfrac12(\mathscr{K}_i+I_j),}
where $\mathscr{K}_i:=I_i\setminus\bigcup\nolimits_n U_n$ and $j\in\{1,2\}$ with $i\neq j$.
\end{enumerate}}

\begin{proof}
First we are going to prove the necessity of cases (i), (ii), and (iii). To do this, assume 
that $(\varphi,f_1,f_2)$ solves equation \eq{E0}.

If $\mathscr{Z}_\varphi$ has an empty interior then, in view of \thm{T0}, $(\varphi,f_1,f_2)$ 
must be extremal, so we obtain solution (i). If $\mathscr{Z}_\varphi$ has a nonempty interior, 
then either $\mathscr{Z}_\varphi=D$ or $\mathscr{Z}_\varphi\neq D$. If $\mathscr{Z}_\varphi=D$, 
then again we get that $(\varphi,f_1,f_2)$ is extremal. Hence, in the sequel, we suppose that 
$\mathscr{Z}_\varphi\neq D$ is valid. Within this sub-case, two main cases can be distinguished: 
either the functions $f_1$ and $f_2$ are both constant on their domain or at least one of them 
is not constant. Since $\mathscr{Z}_\varphi$ is different from $D$, in the first case we get 
that the functions $f_1$ and $f_2$ must take the same value, so we again get case (i). 
Therefore, suppose also that at least one of the functions $f_1$ and $f_2$ is not constant on 
its domain. This means that either (C1) or (C2) is satisfied.

\emph{Case 1. Assume first that neither $f_1$ nor $f_2$ is constant on its domain.} Since 
$\mathscr{Z}_\varphi^c$ is nonempty, the intersection $f_1(I_1)\cap f_2(I_2)$ is nonempty. Then, 
by \rem{R0}, there exist $\lambda,\mu\in f_1(I_1)\cap f_2(I_2)$ such that at least one of the 
inequalities
\Eq{IQ}{
 -\infty<\min(a_1(\lambda),a_2(\lambda))\qquad\text{or}\qquad
		\max(b_1(\mu),b_2(\mu))<+\infty
}
is met. Define the open subintervals
\Eq{iv}{
U_1:=\,]\inf I_1,a_1(\lambda)[,\quad U_2:=\,]b_1(\mu),\sup I_1[,\quad
V_1:=\,]\inf I_2,a_2(\lambda)[,\quad\text{and}\quad V_2:=\,]b_2(\mu),\sup I_2[\,.
}
In view of the inequalities under \eq{IQ}, at least one of the conditions 
$\emptyset\notin\{U_1,V_1\}$ or $\emptyset\notin\{U_2,V_2\}$ is satisfied. Furthermore, by our 
condition (C1), we have $a_i(\lambda)<\sup I_i$, $\inf I_i<b_i(\mu)$, and $a_i(\lambda)\leq 
b_i(\mu)$ hold for any $i\in\{1,2\}$. Consequently, $K_1$ and $K_2$ are indeed proper closed 
subintervals.

For brevity, denote $A:=\frac12(K_1+I_2)\cup\frac12(I_1+K_2)$. The intervals $K_1$ and $K_2$ 
contain at least one point, hence $\emptyset\neq\tfrac12(K_1+K_2)\subseteq A$, which shows that 
$A$ is a nonempty open subinterval of $D$. Now all we need to prove is that the inclusion 
$A\subseteq\mathscr{Z}_\varphi$ is fulfilled. We prove this indirectly.
	
Suppose that we have a point $p\in A$ for which $\varphi(p)\neq 0$. Since $p\in A$, there exist 
$x\in I_1$ and $y\in I_2$ such that $p=\frac12(x+y)$ and that at least one of the inclusions 
$x\in K_1$ or $y\in K_2$ holds. On the other hand, in view of our condition (C1) and of the 
maximality of the intervals under \eq{iv}, we must have $(I_2 \hspace{.3mm}|\hspace{.3mm} 
p)_{I_1}\subseteq I_1\setminus K_1$ and $(I_1 \hspace{.3mm}|\hspace{.3mm} p)_{I_2}\subseteq 
I_2\setminus K_2$ simultaneously. These inclusions yield that $x\in I_1\setminus K_1$ and $y\in 
I_2\setminus K_2$. The contradiction shows that $A\subseteq\mathscr{Z}_\varphi$ is valid. Thus 
we have obtained solution (ii).

\emph{Case 2. Assume now that exactly one of the functions $f_1$ and $f_2$ is constant on its 
domain.} Without loss of generality, we can assume that $f_2(y)=\lambda$ for some 
$\lambda\in\R$ whenever $y\in I_2$. Since the set $\mathscr{Z}_\varphi^c$ is not empty there 
exists a nonempty open subinterval of $I_1$ on which $f_1$ takes the value $\lambda$. By our 
assumption, this interval must be different from $I_1$. This forces that $I_2$ is bounded from 
below or bounded from above.

Let $\mathscr{U}$ consist of the open subintervals of $I_1$ fulfilling the condition that each 
of them is maximal with respect to the property that $f_1$ is identically $\lambda$ on it. On 
the one hand, the family $\mathscr{U}$ is nonempty. On the other hand, having a non-extremal 
solution, $I_1$ does not belong to $\mathscr{U}$. Denote $d:=\diam\,I_2$ and consider the 
subfamily
\Eq{*}{
\mathscr{U}_d:=\{U\in\mathscr{U}\mid d\leq \diam\,U\}.}

Let us first consider the case where the diameter of each member of $\mathscr{U}$ is less than 
$d$, that is, where the subfamily $\mathscr{U}_d$ is empty. Then, for any 
$p\in\cl\,\mathscr{Z}_\varphi^c$, it follows that either
\Eq{Ineq}{\inf(2p-I_i)=2p-\sup I_i<\inf I_j\qquad\text{or}\qquad
\sup I_j<\sup(2p-I_i)=2p-\inf I_i} holds, since $\diam(2p-I_i)=d$. Consequently at least 
one of the disjoint intervals $U_1:=\,]\inf I_j,a_j(\lambda)[$ and $U_2:=\,]b_j(\lambda),\sup 
I_j[$ is nonempty. Moreover, by our assumption, we have 
$\max(\diam\,U_1,\diam\,U_2)<d\leq+\infty$, which yields that $I_j$ is bounded from below or 
bounded from above. In this subcase, define $(U_n):=\{U_1,U_2\}\setminus\{\emptyset\}$. 

If the family $\mathscr{U}_d$ is nonempty, then it can be at most countable. Then let 
$(U_n):=\mathscr{U}_d$.

Now, in any of the above subcases, define $\mathscr{K}_1:=I_1\setminus\bigcup_n U_n$. We show 
that $\varphi$ vanishes on the arithmetic mean $\frac12(\mathscr{K}_1+I_2)$.

Indirectly, assume that this is not the case. Then there exists $p\in\frac12(\mathscr{K}_1+I_2)$ 
such that $\varphi(p)\neq 0$. Then, on one hand, $f_1$ takes $\lambda$ on $(I_2 
\hspace{.3mm}|\hspace{.3mm} p)_{I_1}$. On the other hand $2p-I_2\subseteq\mathscr{K}_1\subseteq 
I_1$, consequently $(I_2 \hspace{.3mm}|\hspace{.3mm} p)_{I_1}=2p-I_2$ and hence $\diam(I_2 
\hspace{.3mm}|\hspace{.3mm} p)_{I_1}=d$. In the case $\mathscr{U}_d=\emptyset$ this is 
impossible, thus $\mathscr{U}_d$ cannot be empty. Then there exists a member $U_{n_0}$ of the 
sequence $(U_n)$ such that $(I_2 \hspace{.3mm}|\hspace{.3mm} p)_{I_1}\subseteq U_{n_0}$. Hence, 
by the definition of $\mathscr{K}_1$, the interval $(I_2 \hspace{.3mm}|\hspace{.3mm} p)_{I_1}$ 
cannot be a subset of $\mathscr{K}_1$, which is a contradiction. To summarize, in this third 
case we have obtained solution (iii).

Finally we turn to the sufficiency of (i), (ii), and (iii). In fact, the sufficiency of (i) is 
obvious. 

To prove the sufficiency of (ii), let $x\in I_1$, $y\in I_2$, and set $u:=\frac12(x+y)$. If 
$x\in U_1$ and $y\in V_1$ or $x\in U_2$ and $y\in V_2$, then $f_1(x)=f_2(y)$, and hence 
equation \eq{E0} is trivially satisfied. Thus we can assume that this is not the case. Then we 
can distinguish two subcases, namely $x\notin U_1$ or $y\notin V_1$.

If $x\in I_1\setminus(U_1\cup U_2)=K_1$, then $u\in A$ follows, which means that $\varphi(u)=0$, 
hence equation \eq{E0} is satisfied. If $x\notin U_1$ and $y\notin V_2$, then $\inf K_1=\sup 
U_1\leq x$ and $y\leq\sup K_2=\inf V_2$ hold. Observe that
\Eq{*}{
\inf A=\min\big(\inf\tfrac12(K_1+I_2),\inf\tfrac12(I_1+K_2)\big)\,\,\text{ and }\,\,
\sup A=\max\big(\sup\tfrac12(K_1+I_2),\sup\tfrac12(I_1+K_2)\big).}
Using these, we obtain that
\Eq{*}{
\inf A
\leq\tfrac12(\inf K_1+\inf I_2)
\leq u
\leq \tfrac12(\sup I_1+\sup K_2)
\leq\sup A.
}
Since $A$ is an interval, it follows that $u\in A$, and hence $\varphi(u)=0$. Case $y\notin V_1$ 
can be treated in a similar way.

Finally, to prove the sufficiency of (iii), for brevity assume that $i=1$ and let $x\in I_1$, 
$y\in I_2$, and $u:=\frac12(x+y)$. If $x\in U_n$ for some index $n$, then 
$f_1(x)=\lambda=f_2(y)$, thus equation \eq{E0} is satisfied. Otherwise $x\in\mathscr{K}_1$, 
which yields that $\varphi(u)=0$.
\end{proof}

As a consequence of \thm{M}, we obtain Theorem 3 of paper \cite{KisPal19}. For the sake of 
clarity, the statement in question is formulated using the terminology of the current paper.

\Cor{*}{
Let $I\subseteq\R$ be a nonempty open subinterval, $\varphi:I\to\R$ be such that 
$\mathscr{Z}_\varphi$ is closed in $I$, and $f:I\to\R$. Then the pair $(\varphi,f)$ is a 
solution of equation
\Eq{old}{
\varphi\Big(\frac{x+y}2\Big)\big(f(x)-f(y)\big)=0,\qquad x,y\in I
}
if and only if either $f$ is constant on $I$ and $\varphi$ is any function or there exist 
constants $\lambda,\mu\in\R$ and open subintervals $U,V\subseteq I$ with
\Eq{*}{
\inf U=\inf I\quad\text{or}\quad
\sup V=\sup I
}
such that
\Eq{oldi}{
f(x)=\lambda\quad\text{if}\quad x\in U,\qquad
f(y)=\mu\quad\text{if}\quad y\in V,\qquad\text{and}\qquad
\varphi(u)=0\quad\text{if}\quad u\in\tfrac12(I+I\setminus(U\cup V)).
}}

\begin{proof} Apply \thm{M} for $I:=I_1:=I_2$ and $f:=f_1:=f_2$.
\end{proof}

%\bibliography{publ,funcequ,plus}
%\bibliographystyle{plain}

\end{document}